\documentclass[10pt,letterpaper]{article}
\usepackage[top=0.85in,left=2.75in,footskip=0.75in]{geometry}

\usepackage{amsmath,amssymb}

\usepackage{changepage}

\usepackage{textcomp,marvosym}

\usepackage{cite}

\usepackage{nameref,hyperref}

\usepackage[right]{lineno}

\usepackage{microtype}
\DisableLigatures[f]{encoding = *, family = * }

\usepackage[table]{xcolor}

\usepackage{array}

\usepackage{caption}
\usepackage{subcaption}
\usepackage[textsize=tiny]{todonotes}
\usepackage{nicefrac}
\usepackage{siunitx}

\newcommand{\EFF}{{E\hspace{-.15em}F\hspace{-.1em}F}}

\newcommand{\eff}{\phi_{\EFF}}

\newcommand{\pro}{\textup{pro}}

\newcolumntype{+}{!{\vrule width 2pt}}

\newlength\savedwidth

\newcommand\thickhline{\noalign{\global\savedwidth\arrayrulewidth\global\arrayrulewidth 2pt}%
\hline
\noalign{\global\arrayrulewidth\savedwidth}}

\raggedright
\usepackage[aboveskip=1pt,labelfont=bf,labelsep=period,justification=raggedright,singlelinecheck=off]{caption}

\makeatletter
\renewcommand{\@biblabel}[1]{\quad#1.}
\makeatother

\usepackage{lastpage,fancyhdr,graphicx}
\usepackage{epstopdf}
\fancyheadoffset[L]{2.25in}
\fancyfootoffset[L]{2.25in}
\begin{document}
\vspace*{0.2in}

\begin{flushleft}
{\Large
\textbf\newline{A phase-field model for non-small cell lung cancer under the effects of immunotherapy} %
}
\newline
\\
Marvin Fritz\textsuperscript{1},
Christina Kuttler\textsuperscript{2},
J. Tinsley Oden\textsuperscript{3},
Pirmin Schlicke\textsuperscript{2},
Christian Schumann\textsuperscript{4},
Andreas Wagner\textsuperscript{2,*},
Barbara Wohlmuth\textsuperscript{2}

\bigskip
\textbf{1} Computational Methods for PDEs, Johann Radon Institute for Computational and Applied Mathematics, Linz, Upper Austria, Austria,
\textbf{2} School of Computation, Information and Technology, Technical University of Munich, Munich, Bavaria, Germany,
\textbf{3} Oden Institute for Computational Engineering and Sciences, The University of Texas at Austin, Austin, Texas, United States of America,
\textbf{4} Clinic of Pneumology, Thoracic Oncology, Sleep and Respiratory Critical Care, Klinikverbund Allg\"au, Kempten, Bavaria, Germany,
\bigskip

* \url{wagneran@cit.tum.de}

\end{flushleft}
\section*{Abstract}
Formulating tumor models that predict growth under therapy is vital for improving patient-specific treatment plans.
In this context, we present our recent work on simulating non-small-scale cell lung cancer (NSCLC) in a simple, deterministic setting for two different patients receiving an immunotherapeutic treatment.

At its core, our model consists of a simple Cahn-Hilliard-based phase-field model describing the evolution of proliferative and necrotic tumor cells. These are coupled to a simple nutrient model that drives the growth of the proliferative cells and their decay into necrotic cells. A single scalar value represents the immunotherapeutic agents in the entire lung, which decreases the proliferative cell concentration during therapy. An ordinary differential equation (ODE) model describes their evolution. Finally, reaction terms provide a coupling between all these equations.
By assuming spherical, symmetric tumor growth and constant nutrient inflow, we simplify this full 3D cancer simulation model to a reduced 1D model.

We can then resort to patient data gathered from computed tomography (CT) scans over several years to verify our model. For the reduced 1D model, we show that our model can qualitatively describe observations during immunotherapy by fitting our model parameters to existing patient data. Our model covers cases in which the immunotherapy is successful and limits the tumor size, as well as cases predicting a sudden relapse, leading to exponential tumor growth.

We then move from the reduced model back to the full 3D cancer simulation in the lung tissue. Thereby, we show the predictive benefits a more detailed patient-specific simulation including spatial information could yield in the future.

\section*{Author summary}
Lung cancer is one of the deadliest diseases, with low long-term survival rates.
Its treatment is still very heuristic since patients respond differently to the same treatment plans.
Therefore,  patient-specific models for predicting tumor growth and the treatment response are necessary for clinicians to make informed decisions about the patient's therapy and avoid a trial and error based approach.
We make a small step in that direction by introducing a model for simulating cancer growth and its treatment inside a 3D lung geometry. In this model, we represent tumor cells by a volume fraction field that varies over space and time. We describe their evolution by a system of partial differential equations, which include patient- and treatment-specific parameters capturing the different responses of patients to the therapies. Our simulation results are compared to clinical data and show that we can quantitatively describe the tumor's behavior with some parameter set. This enables us to change therapies and analyze how these changes could have impacted the patient's health.

\section*{Introduction}

A major challenge of mathematical oncology is predicting the growth of tumors \cite{rockne2019introduction}. 
Cancer is a class of diseases characterized by numerous point mutations in the genome that result in the uncontrolled growth and spread of cells. Overviews of its biological details are, e.g., presented in~\cite{weinberg2006biology,graham2017measuring}.
 The body's immune system can suppress tumor growth by inhibiting cell growth or by destroying cancer cells. On the other hand, it can also promote tumor progression by selecting tumor cells that are better able to survive in an immunocompetent host or by establishing conditions within the tumor microenvironment that facilitate tumor outgrowth~\cite{hanahan2011hallmarks,schreiber2011immunoediting}. 
Immunotherapy attempts to boost the body's immune system and immune responses against cancerous cells to eliminate them and understanding this ability has revolutionized the field of oncology~\cite{zhang2020history}. For a variety of cancer types, immunotherapy has been proven to feature a significant clinical benefit in patients with advanced stages of cancer and is, as of today, well-established as a standard treatment method~\cite{rounds2015nivolumab,keating2016nivolumab,iwai2017cancer}. 
However, it is extremely difficult to accurately identify spatial tumor growth on a patient-specific level, especially under a treatment plan~\cite{ghaffari2022classical,ezhov2023learn}. 
A considerable variety of mathematical models has helped to improve the understanding of biological principles in cancer growth and treatment response \cite{laird_dynamics_1965,norton_gompertzian_1988,benzekry_classical_2014,bilous_quantitative_2019,schlicke2021mathematical} and their predictive power \cite{benzekry_development_2021,Benzekry2023mechanistic}.
Globally, lung cancer is the leading cancer-related mortality cause. Roughly 85\% of all lung cancer cases are non small cell lung cancer (NSCLC). Its five year survival probability is approximately 22\%~\cite{bray2018global}.\\

For an overview of how to model spatial tumor growth, including continuum, discrete, and hybrid models, we refer to the reviews in \cite{lowengrub2009nonlinear,cristini2010multiscale}.
One class of continuum models relies on simple reaction-diffusion-equations and elastic models for the mass effect~\cite{clatz2005realistic,subramanian2019simulation,bowers2020characterization}. 
We will rely on phase-field models commonly used for modeling cell dynamics, since they allow us to model the observed cell-to-cell adhesion \cite{friedl2009collective} between tumor cells by energy terms.
Simple two-species models consisting of a phase-field model for the tumor and a nutrient equation of reaction-diffusion type are introduced and analyzed in \cite{garcke2016cahn,frigeri2015diffuse}.
Models, including increasingly more complicated flow models, are given in \cite{ebenbeck2019analysis,ebenbeck2019cahn,fritz2019unsteady,lam2017thermodynamically,lorenzo2022identifying}.
A more general theoretical approach to multispecies models can be found in \cite{garcke2018multiphase,oden2010general,wise2008three,cristini2009nonlinear,frieboes2010three}.
The number of modeled species and fields strongly depends on the choice of the particular studied effects.
More specialized large models can be found in \cite{lima2014hybrid,hawkins2012numerical,fritz2021modeling,fritz2021analysis,lorenzo2016tissue}.

Regarding spatial models, including cancer therapy \cite{powathil2013towards} introduces a hybrid model to study the impact of different chemo- and radiation therapies on tumor cells.
The chemotherapy of breast cancer with a drug-delivery model is discussed in \cite{wu2022towards}, and with a reaction-diffusion model in \cite{lorenzo2022identifying,jarrett2018incorporating}.
In \cite{rockne2015patient}, a continuum model specialized in radiation therapy for brain tumors is discussed. 
For prostate cancer, a combination of chemotherapy with antiangiogenic therapy and the optimization of the treatment are given in \cite{colli2021optimal}.
Chemotherapy is also included in the multispecies phase-field model approach of \cite{fritz2021subdiffusive}.\\
\ \\
The present work aims to develop a mathematical model for the spatial growth behavior of solid tumors in NSCLC in a simple, deterministic setting.
It addresses the effects of immunotherapy application and allows the description of its influences on the tumor's spatial structure.
The model framework is applied to two data sets acquired from clinical patients that have shown qualitatively different therapy outcomes. Analysis of the model sheds light on the corresponding parameter relations that determine different clinical outcomes that could potentially be estimated in a prognostic setting to improve the prediction of clinical outcomes and therapy choices.

Here, passive immunotherapy is considered that uses monoclonal antibodies to improve the immune response by regulating T-cell activity in the effector phase of the immune response via the so-called PD-1/PD-L1-pathway. The downregulation usually caused by PD-1 activation prevents collateral damage during an immune response~\cite{quezada2013exploiting,ribas2012tumor}. However, tumor cells can impair this pathway by expressing the corresponding ligands PD-L1 and PD-L2 that bind to the T-cells PD-1 receptor and inactivate the T-cell to decrease the immune response towards the tumor cells~\cite{pardoll2012blockade,quezada2013exploiting,rozali2012programmed,patel2015pd-l1}. Cells expressing the mentioned ligands are targets for, e.g., the drugs Nivolumab and Pembrolizumab with which patients of our clinical data set were treated. 

\section*{Model}\label{sec:model}
Our tumor model consists of the volume fraction of tumor cells $\phi_T$, which we divide into the two cell species of proliferative tumor cells $\phi_P$, and necrotic tumor cells $\phi_N$, such that $\phi_T = \phi_P + \phi_N$.
To keep the model as simple as possible, we do not explicitly model hypoxic cell species.
Furthermore, we introduce the nutrient concentrations $\phi_{\sigma,v}$ and $\phi_{\sigma,i}$ inside the vasculature and interstitial space. 
The concentration of immunotherapeutic agents is given by $\phi_\tau$.
We will first start with the full 3D-Model and then introduce its simplifications for a 1D-Model.

\paragraph*{The full 3D-Model}
We use a generalized Cahn--Hilliard model as in~\cite{lima2014hybrid,fritz2021analysis} with reaction terms to describe the evolution of proliferative cells. By introducing the chemical potential $\mu_P$, the model is characterized by the  fourth-order equation given as 
\begin{equation}
    \begin{aligned}
		\partial_t \phi_P ={}& \nabla\cdot(c_m m_P(\phi_P, \phi_T) \nabla \mu_P) + S_P(\phi_P, \phi_{\sigma,i}, \phi_\tau) 
  \quad\textmd{ in } \Omega, \\
		\mu_P ={}& \partial_{\phi_P} \Psi(\phi_P, \phi_T) - \varepsilon_P^2 \Delta \phi_P - \chi \phi_{\sigma,i}
  \quad\textmd{ in } \Omega,
    \end{aligned}
    \label{eq:model_phi_P}
\end{equation}
with boundary conditions
\[
    \frac{\partial \phi_P}{\partial n}
    =
    0,
    \quad \textmd{ and } \quad
    \frac{\partial \mu_P}{\partial n}
    =
    0
    \quad \textmd{ on } \partial\Omega,
\]
where  $\Psi(\phi_P, \phi_T) = \tilde{\Psi}(\phi_P) + \tilde{\Psi}(\phi_T)$, with $ \tilde{\Psi}(\phi) = c_\Psi \phi^2 (1-\phi)^2$ is a double-well potential with scaling factor $c_\Psi$, and $m_P(\phi_P, \phi_T) = \phi_P^2 (1-\phi_T)^2$ with the constant $c_m > 0 $ describing the cell mobility of the tumor cells. 
The reaction term, $S_P$, is given by
\begin{equation}
\begin{aligned}
	S_P(\phi_P, \phi_{\sigma,i}, \phi_\tau) =& \lambda_P^{pro} \phi_{\sigma,i} \left(\phi_P\right)^{\lambda} \ln\left(\frac{1 + \epsilon_{g}}{\phi_P + \epsilon_{g}}\right)
        - \left( \lambda^{\textup{eff}}_\tau \frac{1}{\omega} \right) \phi_P \frac{\phi_\tau}{\phi_\tau + \phi_\tau^{50}}
        \\
        {}&
        - \lambda_{PN}  \mathcal{H}(\sigma_{PN} - \phi_{\sigma,i}) \phi_P
        .
\end{aligned}
\label{eq:S_P}
\end{equation}

The first term on the right-hand side of \eqref{eq:S_P}, is a Gompertzian growth term \cite{gompertz_nature_1825} acting on the tumor interface, with a small parameter $\epsilon_{g}$ for regularization.
We assume that the growth is proportional to the nutrient concentration and scale it by the parameter $\lambda^{pro}_P$, which has to be calibrated beforehand. For increasing $\lambda \geq 1$, the growth of smaller tumor cell concentrations is penalized, which inhibits the tumor from spreading over the entire lung. 

The second term on the right-hand side of \eqref{eq:S_P} models the decay of tumor cells due to immunotherapy and acts on the tumor volume.
Here, $\frac{1}{\omega}$ is the ratio of the lung weight to the body weight, 
$\lambda^{\textup{eff}}_\tau$ describes the patient-specific effect of the immunotherapy on the tumor, and $\phi^{50}_\tau$ is the drug concentration for the half-maximal response.

The last term of \eqref{eq:S_P} models a decay of proliferative cells due to low nutrient concentrations with rate $\lambda_{PN}$. It is activated by the Heaviside function, $\mathcal{H}$, if the nutrient concentration becomes smaller than $\sigma_{PN}$.

Finally, we have a nutrient-dependent term in the equation for $\mu_P$, which describes the chemotactic effect, reflecting that tumor cell growth follows nutrient gradients.

The necrotic cell field is described by a simple spatial ODE model:
\begin{equation}
    \partial_t \phi_N ={} S_N(\phi_{\sigma,i}, \phi_P)
    \quad \textmd{with} \quad 
    S_N(\phi_{\sigma,i}, \phi_P) = \lambda_{PN}  \mathcal{H}(\sigma_{PN} - \phi_{\sigma,i}) \phi_P,
    \label{eq:necrotic}
\end{equation}
where the growth term on the right-hand side mirrors the decay terms of the proliferative cells.

For the nutrients, we follow the approach in \cite{erbert2021coupled} with a stationary diffusion model, and get
\begin{align}
- \nabla \cdot \kappa_v \nabla \phi_{\sigma,v}
&=
- \xi_{va}  \phi_{\sigma,v}
- \eta_{vi}  \phi_{\sigma,v}
+ \eta_{iv}  \phi_{\sigma,i}
+ \beta \delta_\Gamma,
&\textmd{ in } \Omega,
\label{eq:nutrients_v}
\intertext{and}
- \nabla \cdot \kappa_i \nabla \phi_{\sigma,i}
&=
\eta_{vi}  \phi_{\sigma,v}
- \eta_{iv}  \phi_{\sigma,i}
-\alpha_H (1 - \phi_P) \phi_{\sigma,i}
-\alpha_P \phi_P \phi_{\sigma,i}
&\textmd{ in } \Omega,
\label{eq:nutrients_i}
\end{align}
with zero Neumann boundary conditions: $\frac{\partial \phi_{\sigma,v}}{\partial n} = 0$ and $\frac{\partial \phi_{\sigma,i}}{\partial n} = 0$ on $\partial \Omega$.
Static nutrient models are used in \cite{garcke2017well,garcke2017analysis} and are based on the observation that the tumor evolves on a timescale of days and weeks, while the nutrient model works on a much faster time scale.
Here
$- \xi_{va}  \phi_{\sigma,v}$ model the coupling with the arteries, while
$- \eta_{vi}  \phi_{\sigma,v}$ and  $\eta_{iv}  \phi_{\sigma,i}$
model the nutrient exchange between the vasculature and interstitial space.
The source term couples the equation to a 2D manifold $\Gamma$, approximating the surface of arteries and bronchial tubes, by a Delta-Distribution. 
The amount of nutrients entering the system from outside is controlled by the parameter $\beta > 0$.
The nutrient consumption by the healthy cells is modeled by the decay term $-\alpha_H (1 - \phi_P) \phi_{\sigma,i}$, while  $-\alpha_P \phi_P \phi_{\sigma,i}$ models the enhanced consumption by proliferative cells.
Finally, $\kappa_v$ and $\kappa_i$ are the vessel's and interstitial space's scalar diffusion constants.

A simple ODE model governs the immunotherapy concentration:
\begin{equation}
    \begin{aligned}
        \partial_t \phi_\tau ={}&  - \frac{\ln(2)}{t_{1/2}} \phi_\tau + \frac{N_A}{M_\tau} d_\tau \mathbf{1}_{t \in I_\tau},%
    \end{aligned}
    \label{eq:drugs}
\end{equation}
where the first term describes the decay, and the second term is the drug influx due to injections.
The treatment plan is described by the set $I_\tau \subset [0,\infty)$, which contains the times at which the drug is administered. 
Hence, the indicator function acts as an activation function for the source term.
Furthermore, $N_A=6.022140857\cdot 10^{23}$ denotes the Avogadro constant,  $M_\tau$ is the molar mass of the drug, and $d_\tau=0.24$ is the administered dosage of the antibody. %
As done in \cite{schlicke2021mathematical}, we assume that the immunotherapeutic concentration follows the behavior of a Hill--Langmuir equation. 
 All the model parameters are summarized on the left of Table~\ref{tab:growth}.
\begin{table}[ht]
    \begin{adjustwidth}{-2.25in}{0in} 
	\begin{center}
		\caption{{\bf List of model parameters.} }\label{tab:growth}
		\begin{tabular}{|l + l | l + c | c | c | }
			\hline
			&&& \multicolumn{2}{c|}{\textbf{Values 1D}} & \multicolumn{1}{c|}{\textbf{3D}} \\
			\textbf{Parameter} & \textbf{Dimension} & \textbf{Description} 
			& \textbf{\small Patient 1} & \textbf{\small Patient 2} & \textbf{\small Patient 1} \\
			\thickhline
			$c_{m}$ & $\left[ m^2/d \right]$ & Mobility factor of proliferative tumor cells. & \multicolumn{3}{c|}{$5\cdot10^{-4}$} \\
			\hline
		    $\lambda_P^\pro$ & $\left[ 1/d \right]$ & Proliferation rate of proliferative tumor cells. &0.38&0.0038& 0.53\\
			\hline
		    $\lambda$  & $\left[  \right]$ & Power-law for growth of small cell concentrations. &2&1&2\\
			\hline
		    $\lambda_{PN}$ & $\left[ 1/d \right]$ & Decay rate proliferative to necrotic cells. &0.1& 20 & 0.025\\
			\hline
		    $\sigma_{PN}$  & $\left[  \right]$ & Nutrient threshold for necrotic growth. & \multicolumn{3}{c|}{0.2}\\
			\hline
		    $\varepsilon_P$  & $\left[ m \right]$ & Tumor interface width. & \multicolumn{3}{c|}{$5\cdot10^{-4}$}\\
			\hline
		    $c_\Psi$  & $\left[  \right]$ & Scaling double-well potential. & \multicolumn{3}{c|}{2} \\
			\hline
		    $\epsilon_g$  & $\left[ \right]$ & Regularization for Gompertzian growth term. & \multicolumn{3}{c|}{0.1} \\
			\hline
			$\chi$  & $\left[  \right]$ & Chemotaxis factor. & \multicolumn{3}{c|}{0}\\
			\thickhline
			$\kappa_v$  & $\left[ m^2/d \right]$ & Nutrient diffusion factor in the vasculature. &\multicolumn{2}{c|}{-}& $10^{-3}$\\
			\hline
			$\kappa_i$  & $\left[ m^2/d \right]$ & Nutrient diffusion factor in the interstitial space. & \multicolumn{3}{c|}{$10^{-5}$}\\
			\hline
			$\xi_{va}$  & $\left[ 1/d \right]$ & Coupling smaller to larger arteries. &\multicolumn{2}{c|}{-}&0\\
			\hline
			$\eta_{vi}$  & $\left[ 1/d \right]$ & Coupling vasculature to interstitial space. &\multicolumn{3}{c|}{1}\\
			\hline
			$\eta_{iv}$  & $\left[ 1/d \right]$ & Coupling interstitial space to vasculature. &\multicolumn{3}{c|}{0}\\
            \hline
			$\alpha_H$  & $\left[ 1/d \right]$ & Nutrient consumption rate healthy cells. & \multicolumn{3}{c|}{4}\\
            \hline
			$\alpha_P$  & $\left[ 1/d \right]$ & Nutrient consumption rate proliferative cells. & \multicolumn{3}{c|}{1} \\
            \hline
			$\beta$  & $\left[ m/d \right]$ & Source term large vessels. &\multicolumn{2}{c|}{-}& 0.1 \\
			\thickhline
			$\lambda_\tau^\text{eff}$  & $\left[ kg/d \right]$ & Immunotherapeutic effect under application. &4.49& 0.55 & 4.49\\
			\hline
			$\phi_\tau^{50}$ & $\left[ mol \right]$ & Drug concentration for half-maximal response. & \multicolumn{3}{c|}{$1.012\cdot 10^{16}$}\\
			\hline
			$d_\tau$  & $\left[ g \right]$ & Dosage per medication interval. &\multicolumn{1}{c|}{0.24} &\multicolumn{1}{c|}{0.20} &\multicolumn{1}{c|}{0.24}\\
			\hline
			$M_\tau$  & $\left[ g/mol \right]$ & Molar mass of medication. &  \multicolumn{1}{c|}{146.000} & \multicolumn{1}{c|}{143.600} & \multicolumn{1}{c|}{146.000} \\
			\hline
			$t_{1/2}$  & $\left[ d \right]$ & Drug serum half-life time. & \multicolumn{1}{c|}{26.7} & \multicolumn{1}{c|}{22.0} & \multicolumn{1}{c|}{26.7}\\
			\hline
			$\omega$  & $\left[ kg \right]$  & Patient weight. & \multicolumn{3}{c|}{80} \\
			\hline
		\end{tabular}
	\end{center}
	\vspace*{0.25cm}
	\begin{flushleft}
		Left columns: Mathematical symbol. Center columns: Dimension and description. Right columns: Numerical value for our simulation.
        Top rows: Parameters for the tumor growth. Middle rows: Parameters for the nutrient model. Bottom rows: Parameters for the therapy.
	\end{flushleft}
	\end{adjustwidth}
\end{table}

\paragraph*{The reduced 1D-Model}
Since the full 3D model has high computational costs, we assume a spherical symmetric tumor to calibrate our model parameters.
Thus, we begin with a model in which we assume a radial dependency of  $\phi_P = \phi_P(r)$ and $\phi_{\sigma,v} = \phi_{\sigma,v}(r)$, which immediately implies $\mu_P = \mu_P(r)$, $\phi_{\sigma,i}=\phi_{\sigma,i}(r)$ and $\phi_N = \phi_N(r)$.
As a first step, we assume spherical symmetry and, hence, obtain the following simplified model in the radial coordinates:
\begin{align}
    \partial_t \phi_P ={}& 
    \frac{1}{r^2} \frac{\partial}{\partial r} \left( (c_m m_P(\phi_P, \phi_T)) r^2 \frac{\partial \mu_P}{\partial r} \right)
    + S_P(\phi_P, \phi_{\sigma,i}, \phi_\tau) 
      &\quad\textmd{ in } \Omega, \\
            \mu_P ={}& \partial_{\phi_P} \Psi(\phi_P, \phi_T) 
      - \epsilon_P^2 \frac{1}{r^2} \frac{\partial}{\partial r} \left( r^2 \frac{\partial \phi_P}{\partial r} \right)
      - \chi \phi_\sigma
      &\quad\textmd{ in } \Omega,
      \intertext{with boundary conditions}
      \frac{\partial \phi_P}{\partial r}
      &=
      0,
      \quad \textmd{ and } \quad
      \frac{\partial \mu_P}{\partial r}
      =
      0
      &\textmd{ on } \partial\Omega.
\end{align}
The equations for the necrotic cells and immunotherapeutic agents stay unchanged.
For the nutrient model, we assume that the nutrient concentration $\phi_{\sigma,v}$ can be approximated by a pre-given constant value and use 
\begin{align}
- \kappa_i 
\frac{1}{r^2} \frac{\partial}{\partial r} \left( r^2 \frac{\partial \phi_{\sigma,i}}{\partial r} \right)
&=
\eta_{vi}  \phi_{\sigma,v}
- \eta_{iv}  \phi_{\sigma,i}
-\alpha_H (1 - \phi_P) \phi_{\sigma,i}
-\alpha_P \phi_P \phi_{\sigma,i}
&\textmd{ in } \Omega,
\end{align}
for the nutrients in the interstitial space.

We remark that the $r^{-2}$ term does not pose any numerical issues for finite elements since we have to apply the spherical volume measure $4 \pi r^2 \textmd{d}r$ when we bring the system into its weak form and, thus, the $r$-terms cancel.

\section*{Data}\label{sec:patient-data}

Datasets from two clinical patients were available to verify our model. The examined patient data consist of anonymized volumetric measurements of primary tumors and metastases in patients with NSCLC and was acquired from routinely generated CT slices and computations of the secondary appraisal environment \textit{syngo.CT} Lung Computer Aided Detection (LCAD) workflow of Siemens Healthineers, provided in the \textit{syngo.Via} framework.\footnote{The patients were treated in the Clinic of Pneumology, Thoracic Oncology, Sleep and Respiratory Critical Care of the Klinikverbund Allg{\"a}u in Germany. 
The ethics commission of BLAEK (Ethik-Kommission der Bayerischen Landes\"arztekammer), reference number 19021, approved the use of the data.}
As measurements, we have the manually segmented tumor volume, the response evaluation criteria in solid tumors (RECIST v1.1,~\cite{eisenhauer2009new}), and the bidimensional World Health Organization criteria (WHO), all indicating the tumor response with respect to therapy in different defined ways.
\vspace{0.2cm}
\paragraph{Patient 1}
\begin{figure}
    \begin{adjustwidth}{-2.25in}{0in} 
    \includegraphics[width=.5\textwidth+1.125in]{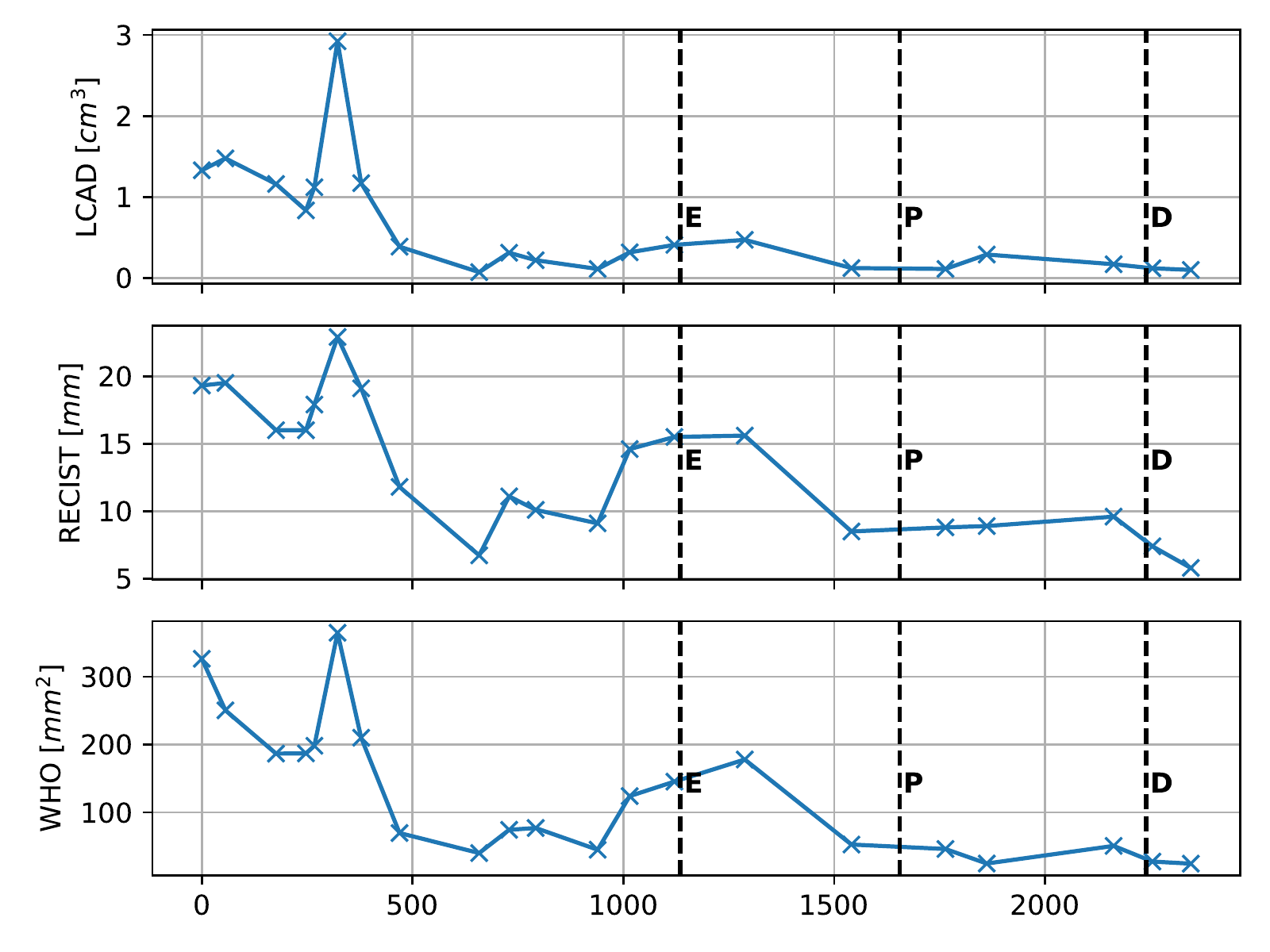}%
    \includegraphics[width=.5\textwidth+1.125in]{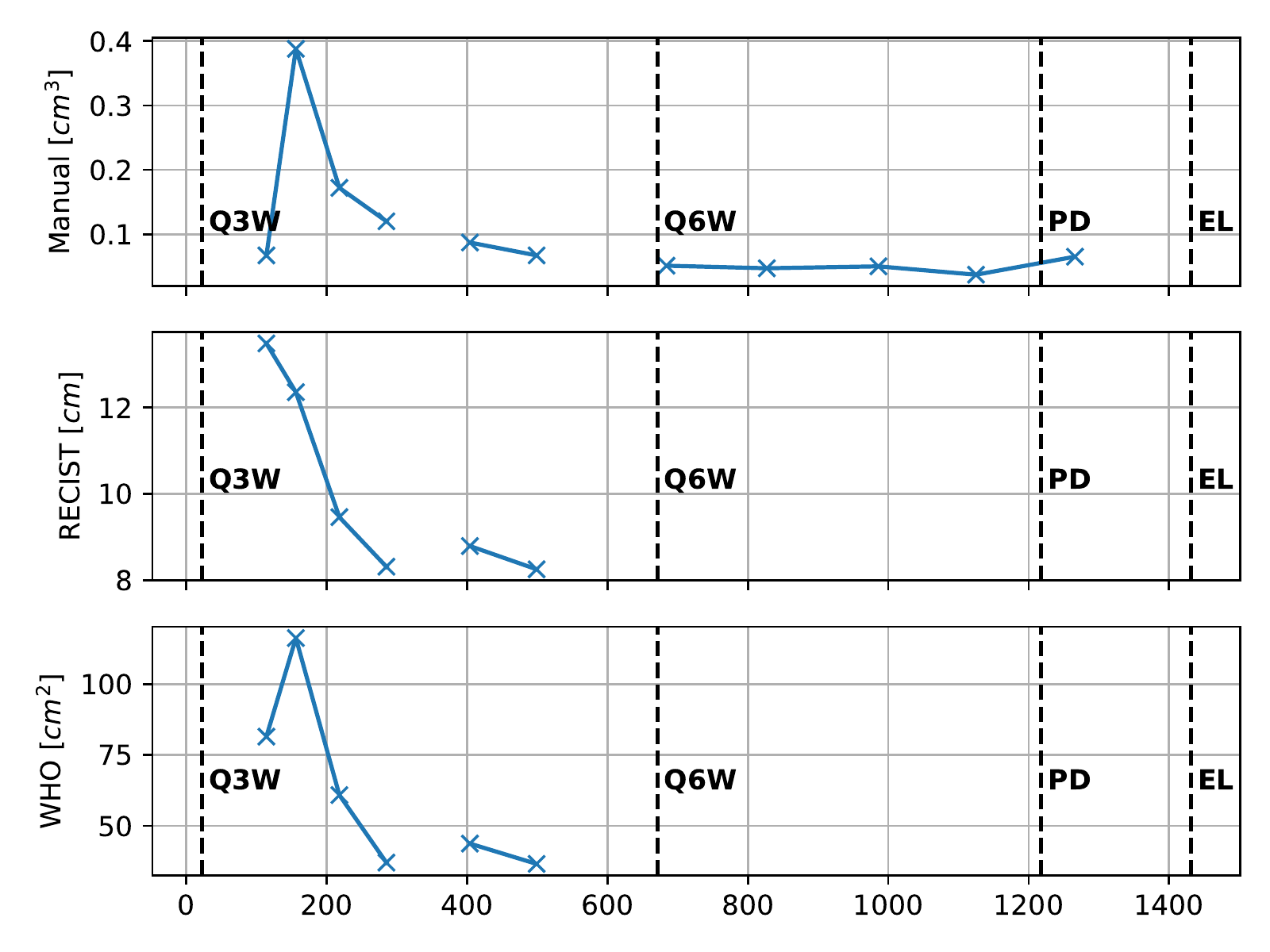}
    \begin{minipage}{.5\textwidth+1.125in}
    \centering (A)
    \end{minipage}%
    \begin{minipage}{.5\textwidth+1.125in}
    \centering (B)
    \end{minipage}
    \parbox{\linewidth}{
    \caption{{\bf Measured patient data from CT scans.}
    Top row: tumor volume. Middle row: RECIST value. Bottom row: WHO value. Vertical dashed lines indicate events that might affect the data quality or therapy. For Patient 1 {\bf E} indicates an Emblie, {\bf P} a therapy break, and {\bf D} the treatment with Dexamethason. For Patient 2 {\bf Q3W} is the beginning of a tree weekly drug administration schedule, {\bf Q6W} the change to a six-weekly administration, {\bf PD} the therapy end, and {\bf EL} the Exitus Letalis.
    }
    \label{fig:measured-values-patients}
    }
    \end{adjustwidth}
\end{figure}
The patient data during the therapy is depicted in Fig.~\ref{fig:measured-values-patients}A, which contains the tumor volume and the RECIST and WHO values.
The tumor was diagnosed at time $t=0$, and the therapy started at time $t_{S}=296$ days.
The patient received Nivolumab as an immunotherapeutic antibody every two weeks, for 30-60 minutes, in an intravenous dose administration from this time on.
At time point (E), the patient developed a temporary lung emboly, significantly influencing the data quality and LCAD accuracy. 
From (P) on, a break of immunotherapeutic drug administration was conducted due to the evident presence of brain metastases.
At that very time point, the treating medical doctors diagnosed a complete remission for the primary tumor in the lung.
All in all, we model our therapy as
\begin{align}
    I_\tau := \left\{ t  :  t_S+14\cdot k \leq t  \leq t_S+14\cdot k + \frac{1}{24}, \textmd{ with } t \leq t_P \right\}.
\end{align}

\paragraph{Patient 2}
For our second patient, the therapy starts before we have our first data point.
The patient receives 200 mg of Pembrolizumab every 3 weeks.
The therapy plan changes at time $t_{Q6W}=514$ days, where the application interval is adjusted from three to six weeks.
Finally, it was terminated at time $t_{PD}=1063$ days with the patient's death (exitus letalis) at time $1306$ days.
The final CT image shows a very large tumor mass that has spread over the entire lung.
The corresponding therapy for this patient is modeled as
\begin{align}
\begin{split}
    I_\tau := &
    \left\{ t  :  t_S+21\cdot k \leq t  \leq t_S+21\cdot k + \frac{1}{24}, \textmd{ with } t \leq t_{Q6W}, k \in \mathbb{N} \right\} \cup
    \\
    &\cup
    \left\{ t  :  t_{Q6W}+42\cdot k \leq t  \leq t_{Q6W}+42\cdot k + \frac{1}{24}, \textmd{ with } t \leq t_{PD}, k \in \mathbb{N} \right\}.
\end{split}
\end{align}

\section*{Results}\label{sec:simulation-results-3d}
In this section, we will discuss our simulation results.
All the used parameter values are summarized on the right of Table~\ref{tab:growth}.
In CT images, only a small part of the tumor cells are visible. To take this into account, we define the visible tumor volume from our simulation as $V_{vis} = \int_\Omega \mathbf{1}_{\{x \in \Omega : \phi(x) > 0.3\}}(x) \textmd{d}x $. Numerically, we apply a time step size of $\nicefrac{1}{24}$ for the tumor and 32 substeps for the evolution of the immunotherapeutic agents.
For the 1D model, we use spatially 500 elements in a spherical domain with radius \SI{4}{\cm}.

\paragraph*{1D results}
\begin{figure}
    \centering
    \begin{adjustwidth}{-2.25in}{0in} 
    \includegraphics[width=.5\textwidth+1.125in]{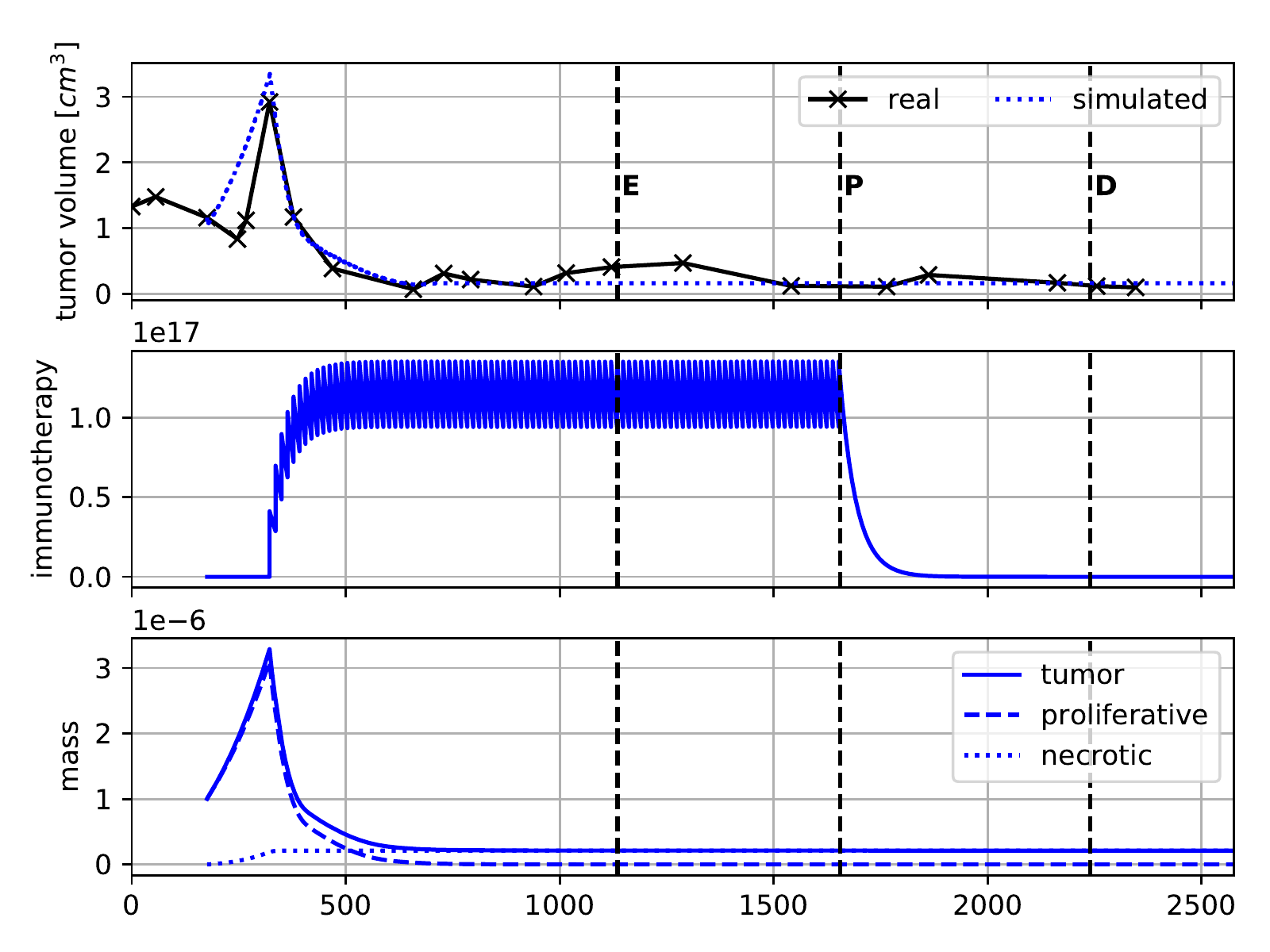}
    \includegraphics[width=.5\textwidth+1.125in]{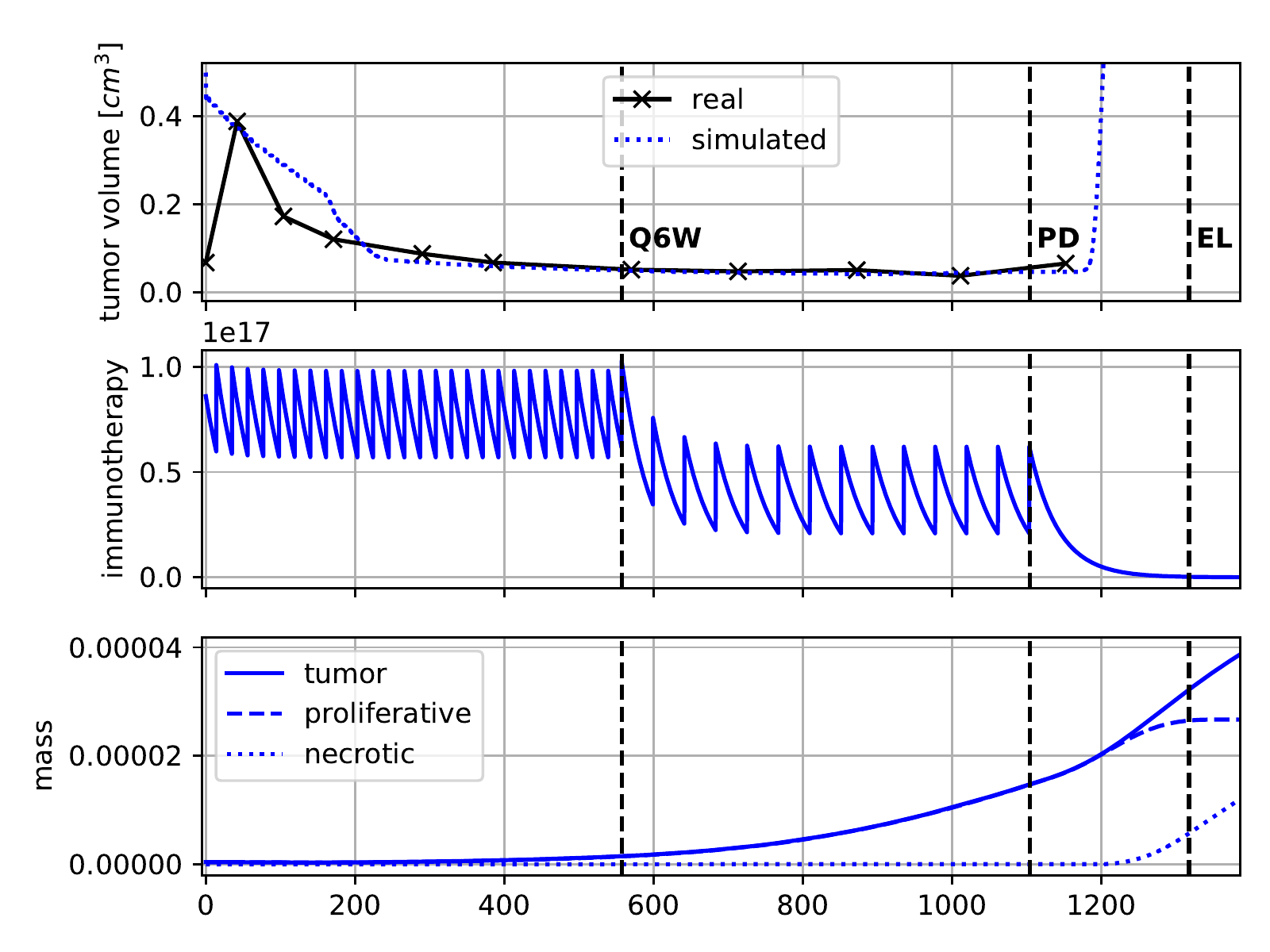}
    \begin{minipage}{.5\textwidth+1.125in}
    \centering (A)
    \end{minipage}%
    \begin{minipage}{.5\textwidth+1.125in}
    \centering (B)
    \end{minipage}
    \parbox{\linewidth}{
    \caption{{\bf Patient data and simulation results:} Contains the visible tumor volume, the number of immunotherapy agents, and the tumor mass over time.
    A: Simulation results for Patient 1. B: Simulation results for Patient 2.
    }
    \label{fig:results_patients}
    }
    \\
    \end{adjustwidth}
\end{figure}
We first discuss the results obtained with the simplified 1D model in spherical coordinates.
Fig.~\ref{fig:results_patients}A shows a comparison of the results obtained for the first patient.
The manually calibrated parameters are $\lambda^{\pro}_P = 0.38, \lambda^{\textmd{eff}}_\tau = 4.49$, $\lambda = 2$ and $\lambda_{PN} = 0.1$.
We start with an initial tumor volume with enough temporal distance from the previous therapy window.
At the top, the visible tumor volume of our simulation is shown as a dotted line, which qualitatively follows the segmented data points.
The number of immunotherapeutic agents is depicted in the middle. Note that after a few months, a quasi-periodic state is reached.
The lower plot depicts the entire tumor mass in the whole domain and therefore includes cell concentrations that cannot be detected yet.
Its evolution mainly depends on the parameter $\lambda$, such that larger values of $\lambda$ impede the tumor cells from spreading over the entire domain.
In this case, the tumor mass curve follows the tumor volume curve.
All in all, the given immunotherapy manages to efficiently eliminate the proliferative tumor cells.
Since, at the end of the therapy only necrotic cells remain, the tumor stays under control, even though no immunotherapeutic agents hamper its growth.

\begin{figure}
    \includegraphics[width=1.\textwidth]{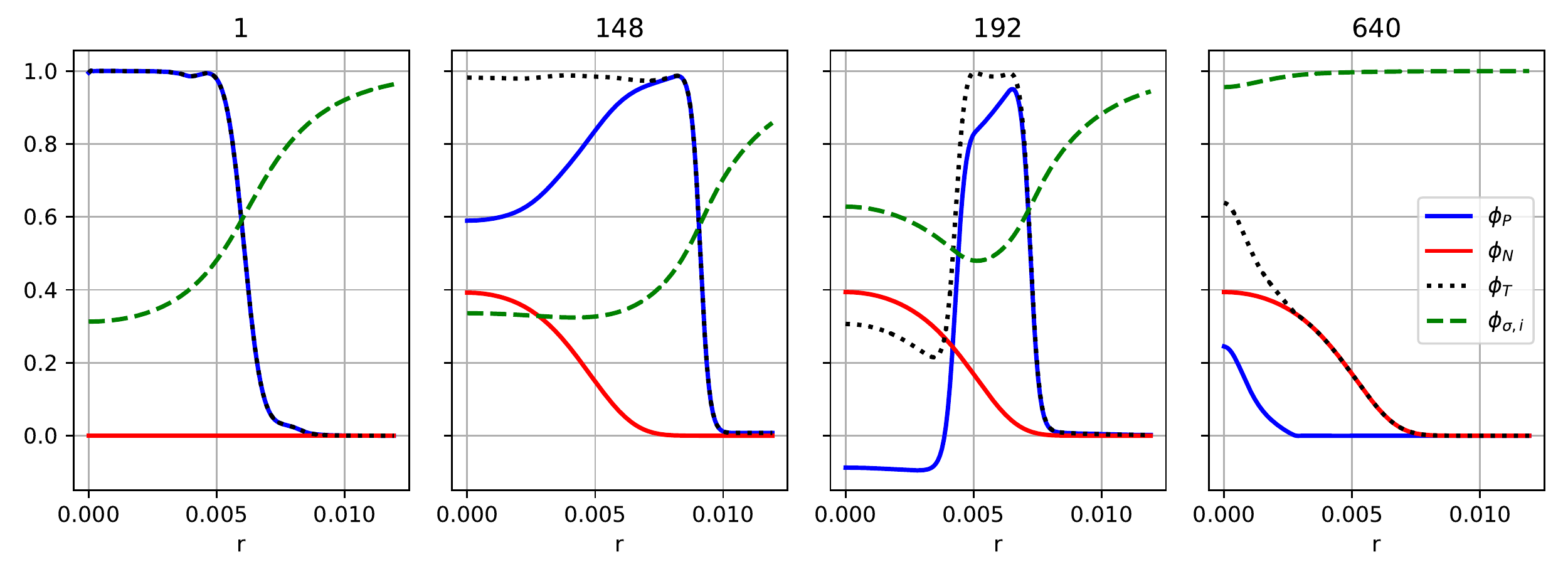}
    \begin{minipage}{.25\textwidth}
    \centering\hspace{0.7cm} (A)
    \end{minipage}%
    \begin{minipage}{.25\textwidth}
    \centering\hspace{0.45cm} (B)
    \end{minipage}%
    \begin{minipage}{.25\textwidth}
    \centering\hspace{0.25cm} (C)
    \end{minipage}%
    \begin{minipage}{.25\textwidth}
    \centering\hspace{0.075cm} (D)
    \end{minipage}
    \parbox{\linewidth}{
    \caption{{\bf Cell- and nutrient fields for Patient 1 plotted over radial distance.}
    A:~Start of the simulation.
    B:~Just before the therapy. 
    C:~Just after therapy starts. 
    D:~During the therapy. 
    }
    \label{fig:results_patient_1_slideshow}
    }
\end{figure}
Fig.~\ref{fig:results_patient_1_slideshow} shows the simulation state in spherical coordinates. In Fig.~\ref{fig:results_patient_1_slideshow}A, the initial state of the tumor simulation is shown. 
The entire tumor consists of proliferative cells, and no necrotic cells have been able to form yet. We see a large drop in the nutrient concentration due to the increased consumption of the tumor cells.
In Fig.~\ref{fig:results_patient_1_slideshow}B, we observe the tumor just before the start of the immunotherapeutic treatment.
Because of the scarcity of nutrients at the tumor core, some of the inner proliferative cells have necrotized.
Since necrotic cells do not deplete nutrients, the nutrient concentrations do not decrease even though the entire tumor had grown.
Fig.~\ref{fig:results_patient_1_slideshow}C is a snapshot taken during the therapy.
Since the immunotherapeutic term acts on the entire volume of proliferative cells, while the growth terms only act on the interfaces, the proliferative cells have increased their interface to the other cell species by becoming active at the tumor boundary, which leads to an unexpected increase of nutrients near the necrotic tumor core.
We also observe a large undershoot of the proliferative cells inside the tumor.
In Fig.~\ref{fig:results_patient_1_slideshow}D, the proliferative tumor shell contracts due to the therapy until the necrotic and proliferative cells form a small tumor sphere. 
Due to numerous nutrients, no conversion of proliferative to necrotic cells happens at this point.
The ongoing immunotherapeutic treatment decreases the proliferative cells until, in the end, only the already present necrotic core remains, which is in agreement with our clinical observations.

Fig.~\ref{fig:results_patients}B shows a comparison of the results obtained for the second patient. 
Here, the therapy is already in progress, and different from the first patient, we have no time window to observe the tumor without any therapy. 
The manually calibrated parameters are $\lambda^{\pro}_P = 0.0038, \lambda^{\textmd{eff}}_\tau = 0.499$, $\lambda = 1$ and $\lambda_{PN} = 20$.
With respect to the tumor volumes in the topmost plot, the first data point is ignored and a plausible tumor volume is assumed, more in line with the previous values.
The visible tumor diminishes in size until it is able to sustain itself with the given nutrients despite ongoing therapy and stays constant. 
The change in the injection schedule from 3 to 6 weeks shows itself in a decreased amount of immunotherapeutic agents, but has no notable impact on the visible tumor.
Even though the visible tumor does not grow, its cells start to spread over the entire lung, which leads to a consistent increase in the total tumor mass in the lower plot.
After the tumor has spread over the entire domain, the visible tumor starts to grow exponentially, which is in line with the observed outcome of the tumor patient. 

\begin{figure}
    \begin{minipage}{1\textwidth}
    \includegraphics[width=1\textwidth]{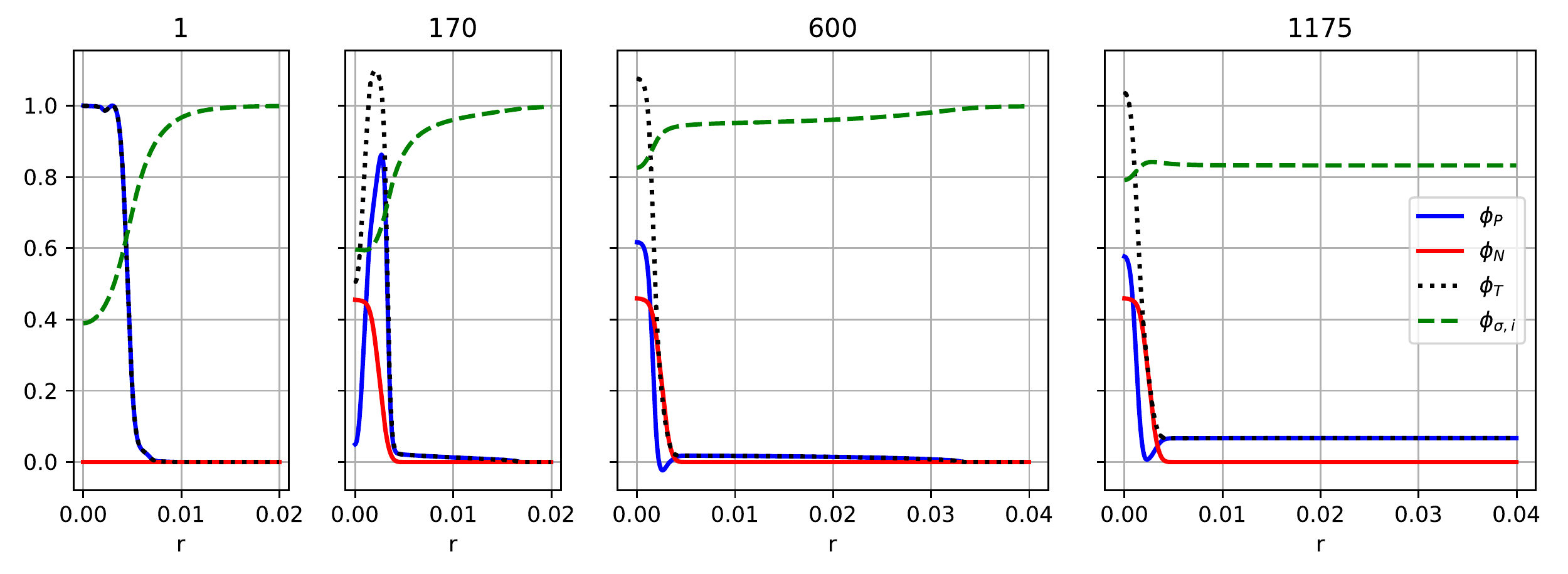}
    \end{minipage}
    \begin{minipage}{.166\textwidth}
    \centering\hspace{0.75cm} (A)
    \end{minipage}%
    \begin{minipage}{.166\textwidth}
    \centering\hspace{0.85cm} (B)
    \end{minipage}%
    \begin{minipage}{.33\textwidth}
    \centering\hspace{0.75cm} (C)
    \end{minipage}%
    \begin{minipage}{.33\textwidth}
    \centering\hspace{0.25cm} (D)
    \end{minipage}%
    \\
    \parbox{\linewidth}{
    \caption{{\bf Cell- and nutrient fields for patient 2 plotted over radial distance.}
    A: Simulation start. B: After 170 days. C: After 600 days. D: After 1.175 days.
    }
    \label{fig:results_patient_2_slideshow}
    }
\end{figure}
Fig.~\ref{fig:results_patient_2_slideshow} gives a more detailed view of what happens spatially.
Fig.~\ref{fig:results_patient_2_slideshow}A describes a similar initial setup to the first patient: The entire tumor consists of proliferative cells. No necrotic core has formed yet, but the nutrient concentration is notably lower inside the tumor.
After 170 days, Fig.~\ref{fig:results_patient_2_slideshow}B is reached, where a necrotic core has formed, the proliferative cells have developed an outer shell, and the nutrient concentration inside the tumor has stabilized.
We also observe that a small amount of proliferative cells has moved away from the primary tumor and spread towards the outer boundary.
In Fig~\ref{fig:results_patient_2_slideshow}C, after 600 days, the proliferative shell has collapsed, and we have a tiny amount of proliferative cells which can sustain themselves due to the increased nutrient concentration.
This configuration is stable over the next year, as we see from Fig.~\ref{fig:results_patient_2_slideshow}D, which is 1.175 days after the start of the simulation. 
This is just before the visible tumor starts to grow again.
We observe that the tumor has spread over the entire domain, which leads to a visible decrease in nutrient concentration from 1 to roughly 0.9.
At this point, the phase separation from the Cahn--Hilliard is initiated, resulting in drastic tumor growth within a short period.
We point out that, especially in the case of explosive growth in the second patient, the simulation geometry does play a role. 
This is the main motivation behind moving from the simplified model in spherical coordinates to more elaborate 3D simulations.

\begin{figure}
    \centering
    \begin{adjustwidth}{-2.25in}{0in} 
    \includegraphics[width=.5\textwidth+1.125in]{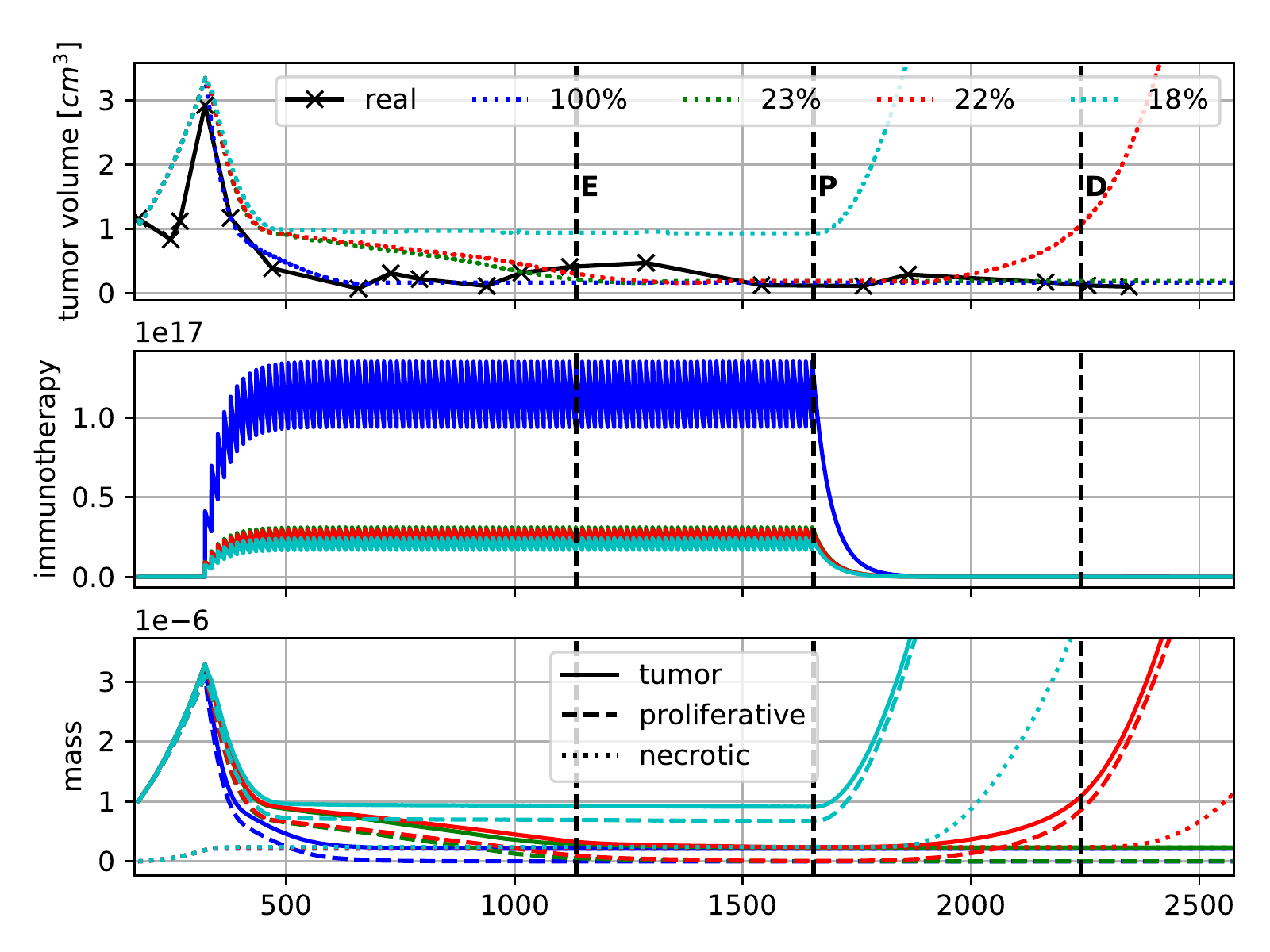}%
    \includegraphics[width=.5\textwidth+1.125in]{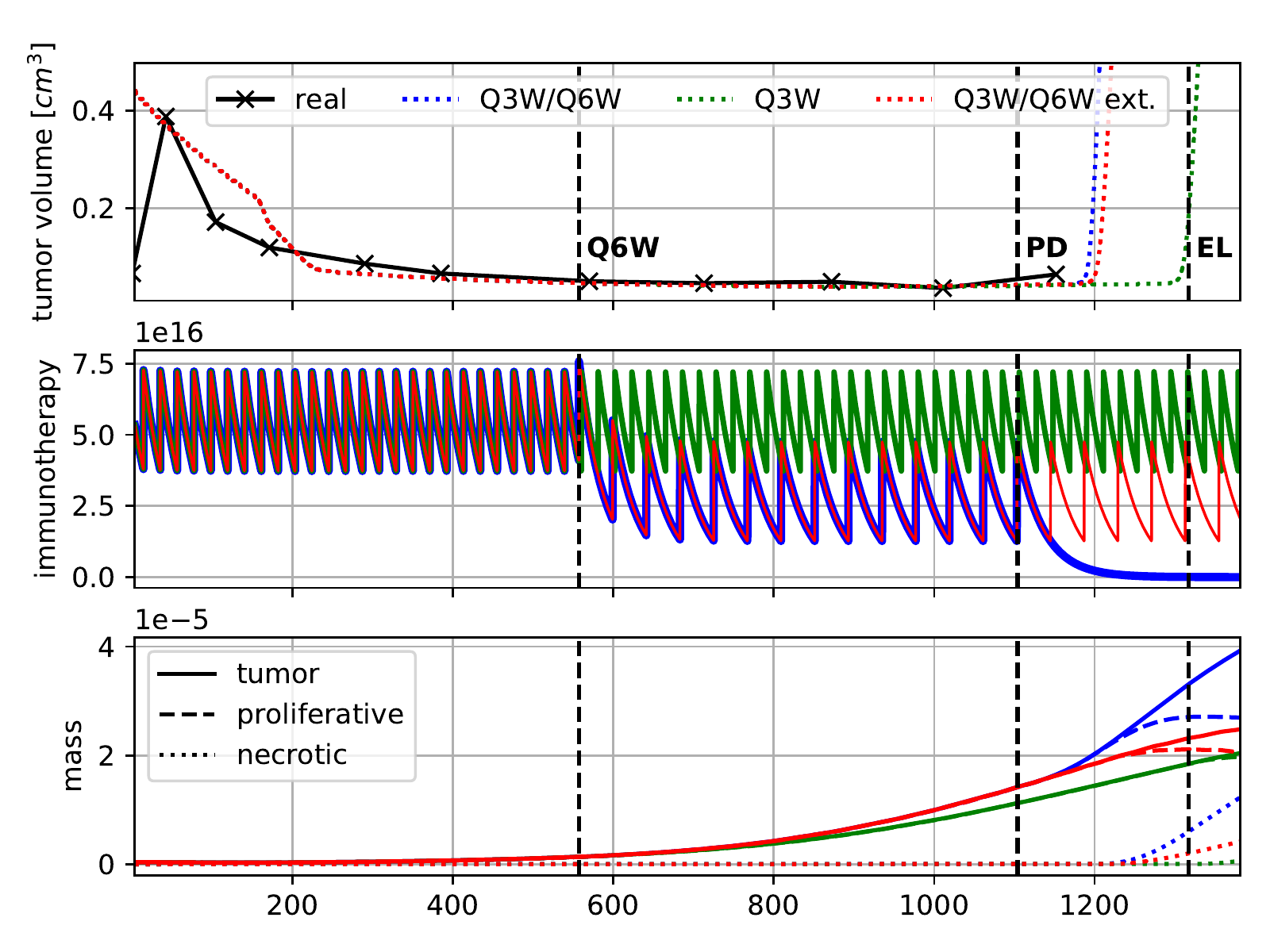}
    \begin{minipage}{.5\textwidth+1.125in}
    \centering (A)
    \end{minipage}%
    \begin{minipage}{.5\textwidth+1.125in}
    \centering (B)
    \end{minipage}
    \\\vspace{0.25cm}
    \parbox{\linewidth}{
    \caption{{\bf Patient data and simulation results:} Contains the visible tumor volume, the number of immunotherapy agents, and the tumor mass over time for different therapies.
    A: Different drug dosages for Patient 1.
    B: Different therapy schedules for Patient 2.
    }
    \label{fig:results_drug_studies}
    }
    \end{adjustwidth}
\end{figure}
Under the assumption that we have meaningful parameter values for both patients, an interesting follow-up question is to investigate how changes in the therapy qualitatively affect its outcomes.

In Fig.~\ref{fig:results_drug_studies}A, we have hypothetically decreased the dosage given to Patient 1 by $77$, $78$, and $82$ percent relative to its original concentration.
At $22$, and $23$ percent of the original dosage, the proliferative cell mass steadily decreases during the therapy. 
While $23$ percent is enough to defeat the tumor before the therapy stops, for $22$ percent, a tiny amount of proliferative cells remain.
These remaining cells then lead to a relapse of tumor growth, such that the tumor exceeds its largest volume even before our simulation time frame is over.
At $18$ percent of the original dosage, we arrive at an equilibrium of the proliferative cells, in which the immunotherapy manages to control its growth but fails to eradicate it.
Thus, as soon as the therapy is terminated, the proliferative tumor starts to grow again at the original rate.
Finally, we note that even considering safety margins, $50$ percent of the dosage would have been more than enough to cure the cancer.

For Patient 2, an interesting question is how the changes in the therapy schedule might have affected the tumor growth.
In Fig.~\ref{fig:results_drug_studies}B, in addition to the original therapy (Q3W/Q6W), we have added the scenario where the patient is given the same dosage continuously every three weeks (Q3W) and the case where the switch to the six weekly cycles was performed, but no therapy break occurred (Q3W/Q6W ext.).
We note that extending the therapy at the given point had close to no effect on the tumor growth.
Even the more aggressive therapy only managed to slow down the growth dynamics by roughly three months before the tumor relapsed.
It appears that according to our model, the given tumor cannot be cured with the used drug, and the therapy only manages to prolong the patient's life by a certain time.
Even significant increases in the dosage are not enough to cure the cancer, since the Hill--Langmuire equation which models the impact of the immunotherapy, approaches a stationary value for $\phi_\tau \to \infty$.
With a continuous three-weekly drug administration scheme, further simulations show that increasing the dosage by a factor of 2 increases the time until the explosive growth occurs by half a year. A factor of 10 could give the patient an additional year.

We will close this part with a short discussion of the parameter values.
Note that the dosage $d_\tau$, the drug's molar mass $M_\tau$ and the serum's half-life time $t_{\nicefrac{1}{2}}$ are known a priori from the given therapy.
The parameters $\lambda^{\pro}_P$, $\lambda$, $\lambda_{PN}$ and $\lambda^{\eff}_\tau$ were assumed patient specific and vary in our examples over a wide range. In particular, $\lambda$ has a strong influence if the tumor stays localized $\lambda=2$ or spreads over the entire domain $\lambda=1$, leading to explosive growth.
The parameters $\lambda^{\pro}_P$ and $\lambda^{\eff}_\tau$ have to counterbalance each other. Especially for Patient 2 where both growth and therapy are always active, they are hard to determine. 

\paragraph*{3D results}
\begin{figure}
    \centering
    \begin{minipage}{\textwidth}
    \includegraphics[width=\textwidth]{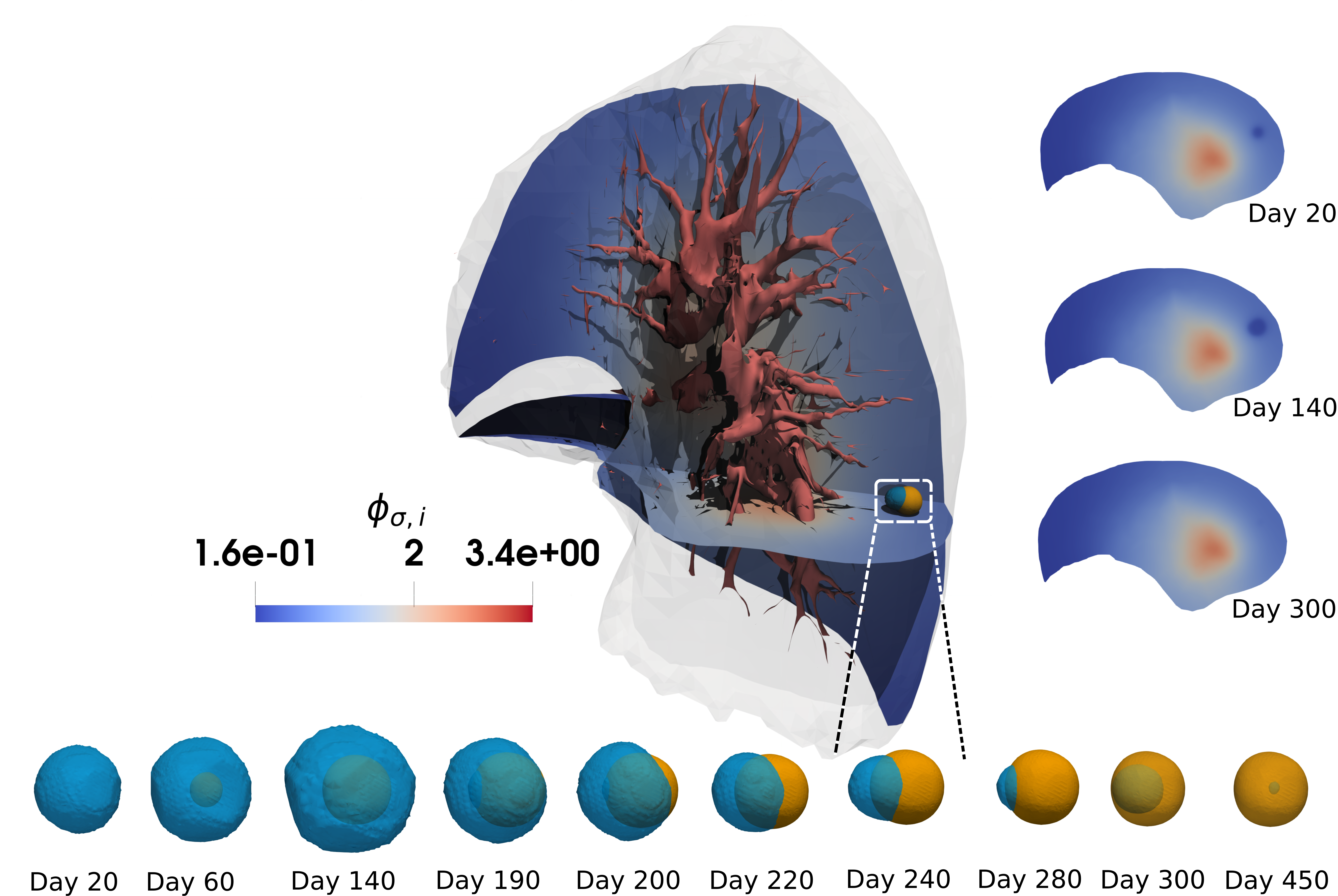}
    \end{minipage}
    \\\vspace*{0.25cm}
    \caption{{\bf 3D tumor simulation:}
    Top-left: Nutrient concentration $\phi_{\sigma,i}$ after 240 days have elapsed. 
    Top-right: Nutrient concentrations at 20, 140 and 300 days.
    Bottom: Tumor shapes consisting of proliferative (blue) and necrotic (orange) cells at different points in time.
    }
    \label{fig:patient_1_3d_tumor_shapes}
\end{figure}
Finally, Fig.~\ref{fig:patient_1_3d_tumor_shapes} shows the tumor evolution for the full 3D model applied to Patient 1.
On the top, a depiction of the lung mask in light gray, of the vasculature $\Gamma$ in red, a vertical and horizontal slice for $\phi_{\sigma,i}$, and the tumor in the lower right are shown.
The proliferative tumor cells are colored blue, while the necrotic cells are depicted in yellow.
A threshold of $0.3$ is chosen for both the contours of the proliferative and necrotic cells.
For the given snapshot in time, the shell of proliferative cells is seen to be oriented in the direction of the nutrient gradient. 
On the right, the nutrient concentration at different times is shown. 
Already at day 20, a nutrient shortage close to the tumor is observed. This area grows with the tumor until it reaches its maximal size after 140 days.
After 300 days, most of the proliferative cells consuming nutrients have disappeared, and the effects of the tumor on the nutrient concentration are no longer visible. 
The tumor evolution for certain points in times are shown.
After 60 days, a small necrotic core has formed. 
It is not at the center of the tumor, but slightly shifted towards the nutrient-poor domain.
The tumor grows until it reaches its maximum size after 140 days, shortly before the immunotherapy starts.
Again, the proliferative cells further away from the nutrients decay until only a tiny proliferative shell remains, which points towards the nutrient source.
Finally, only a few proliferative cells remain inside the tumor, mainly consisting of necrotic cells.
At this point, the visible tumor volume no longer changes and the proliferative core decays further and further.

\section*{Discussion}
Critical microenvironments influencing tumor formation are generally not observable \cite{hanin2018suppression}. 
Deterministic modeling approaches, on the other hand, can aid in understanding the driving dynamics by simulating growth and decline for observable tumor sizes.
It has been demonstrated that our phase-field model can qualitatively describe the tumor volume evolution in the observable window over time for NSCLC patients and can address different outcomes of immunotherapeutic treatment approaches. 
The simulation results of the two clinical cases were explained in a biologically meaningful manner throughout the observation period with results that agree with multiple volumetric measurements at time points acquired during clinical routine examinations. 
This is considered a preliminary proof of concept for the presented model.
In future studies, where the model is calibrated over larger data sets, the parameter estimation should be addressed in particular.
The parameter estimation is then expected to be more stable if prior parameter distributions for general patient cohorts can be identified.\\
We have shown that our model easily generalizes to a full simulation in 3D. 
Backdrawn data on RECIST measurements, which are clinically routine characteristics, is easily processable, but the presented model may contain significantly more clinically relevant information by providing the spatial structure of tumors under therapy influence. \\
The model simulation strongly depends on patient-specific parameters unknown a priori and inferred over the whole simulation time, including the clinical outcomes. 
However, if these parameters can be estimated quickly during treatment or are based on patient-specific clinical covariates, then this approach has prognostic potential.
The long-term therapeutic success and lower but still successful dosages could possibly be determined. 
A preview of this potential is presented in our results, as the qualitative therapeutic outcome can be determined with alternative drug dosage regimens for the two patients presented.\\
This potential should be addressed in future work and be based on large patient data sets. 
Successful early estimation of these parameters allows for an optimal treatment schedule and a minimal drug dosage regimen.
In addition, it could also reduce data acquisition efforts by reducing the number of CT scans during patient treatment to a minimal amount needed for model control, re-calibration, and clinical check-ups.
This has the potential to increase the patient's quality of life on an individual basis. 

\section*{Methods \& Materials}
We will first introduce the numerical discretization of our model in space in time.
Afterwards, the simulation geometry obtained from the CT data is introduced.

\paragraph*{Numerical discretization}

As a spatial discretization, we use finite elements with piecewise-linear basis functions for $\phi_P$, $\mu_P$, $\phi_N$, $\phi_{\sigma,v}$, and $\phi_{\sigma,i}$.
We use a semi-implicit scheme in time for Eq.~\eqref{eq:model_phi_P}.
The source and sink terms, as well as the mobilities, are treated explicitly.
Since the Cahn--Hilliard equation does not guarantee the boundedness, $0 \leq \phi_P \leq 1$, we use a cutoff function on all right-hand side terms to achieve stability.
For the potential $\Psi$, we use a simple convex-concave splitting in time:
Note that the double well potential $\tilde\Psi$ can be decomposed into the $\tilde\Psi_i = \frac{3}{2} c_\Psi \phi^2$ and $\tilde\Psi_e = c_\Psi \left( - 2\phi^3 -\frac{1}{2} \phi^2 + \phi^4 \right)$. 
Clearly, $\Psi_i$ is convex, while $\Psi_e$ stays concave for $\phi \in \left[\frac{1}{2}-\frac{1}{\sqrt{3}}, \frac{1}{2}+\frac{1}{\sqrt{3}} \right] \supset \left[0, 1\right]$.
In our case, this motivates the decomposition $\Psi_i (\phi_P, \phi_N) = \tilde{\Psi}_i(\phi_P) + \tilde{\Psi}_i(\phi_P + \phi_N) $ and $\Psi_e (\phi_P, \phi_N) = \tilde{\Psi}_e(\phi_P) + \tilde{\Psi}_e(\phi_P + \phi_N) $. 
For a fixed $\phi_N$, the latter is concave if $\phi_P \in \left[ \frac{1}{2}-\frac{1}{\sqrt{3}}, \frac{1}{2}+\frac{1}{\sqrt{3}} - \phi_N \right] \supset \left[0, 1 - \phi_N \right]$.
In case of no source terms, we achieve unconditional gradient stability \cite{eyre1998unconditionally} if $\Psi_i$ is treated implicitly and $\Psi_e$ explicitly and the bounds of $\phi_P$ are satisfied.
As a linear solver, we use MINRES with the block-diagonal preconditioner 
\begin{align}
    \mathcal{P} :=
    \begin{pmatrix}
        (c_m \tau) K_m + \sqrt{c_m \tau} M & \\
        & \epsilon^2 K_1 + \frac{6 c_\Psi}{\sqrt{c_m \tau}} M
    \end{pmatrix}
    \label{eq:preconditioner}
    ,
\end{align}
which is a generalization to the preconditioner in \cite{brenner2018robust}. 
For the inversion of the diagonal blocks, we resort to algebraic multigrid~\cite{petsc-user-ref}.

For Eq.~\eqref{eq:necrotic} we use an explicit Euler scheme, while in Eqs.~\eqref{eq:nutrients_v}~and~\eqref{eq:nutrients_i} $\phi_P$ is kept explicit, and we solve for $\phi_{\sigma,v}$ and $\phi_{\sigma,i}$.
Since in our studies, we set $\eta_{iv} = 0$, both equations are decoupled and can be solved separately with a Conjugate Gradient method preconditioned by algebraic multigrid~\cite{petsc-user-ref}. 

Finally, we apply a simple explicit Euler scheme for the ODE model of Eq.~\eqref{eq:drugs}. 
Since Eq.~\eqref{eq:drugs} runs physically on a much smaller time scale, we run it with a much smaller time step size.
The coupling does not pose any difficulties, since Eq.~\eqref{eq:drugs} is completely decoupled and computationally cheap to solve.

For our implementation, we use the FEniCS framework~\cite{logg2012automated}.

\paragraph*{Simulation geometry}
The lung geometry and vasculature structure were extracted from the patient's CT data (see Patient Data) at the therapy's onset.
For the lung geometry, the extraction was performed with \emph{3D Slicer} \cite{Kikinis2014}, which was simplified with \emph{3D Builder} by Microsoft Corporation and \emph{blender} \cite{blender} before remeshing it coarsely with \emph{gmsh} \cite{geuzaine2009gmsh}.
The vasculature was extracted using slicers' \emph{vmtk} extension \cite{vmtk,antiga2008image}.
\begin{figure}
    \begin{adjustwidth}{-2.25in}{0in} 
    \begin{minipage}{.5\textwidth+1.125in}
    \centering
    \includegraphics[width=0.70\textwidth]{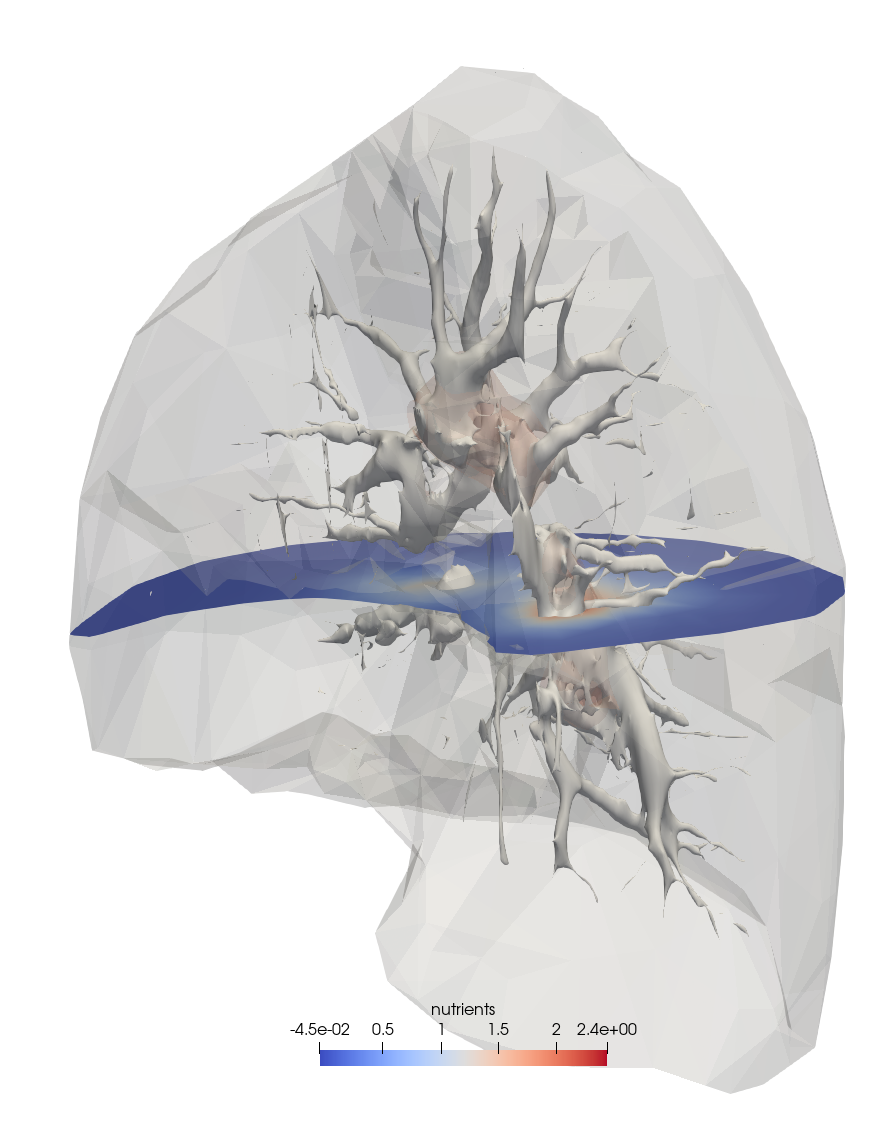}%
    \end{minipage}%
    \begin{minipage}{.5\textwidth+1.125in}
    \centering
    \includegraphics[width=0.56\textwidth]{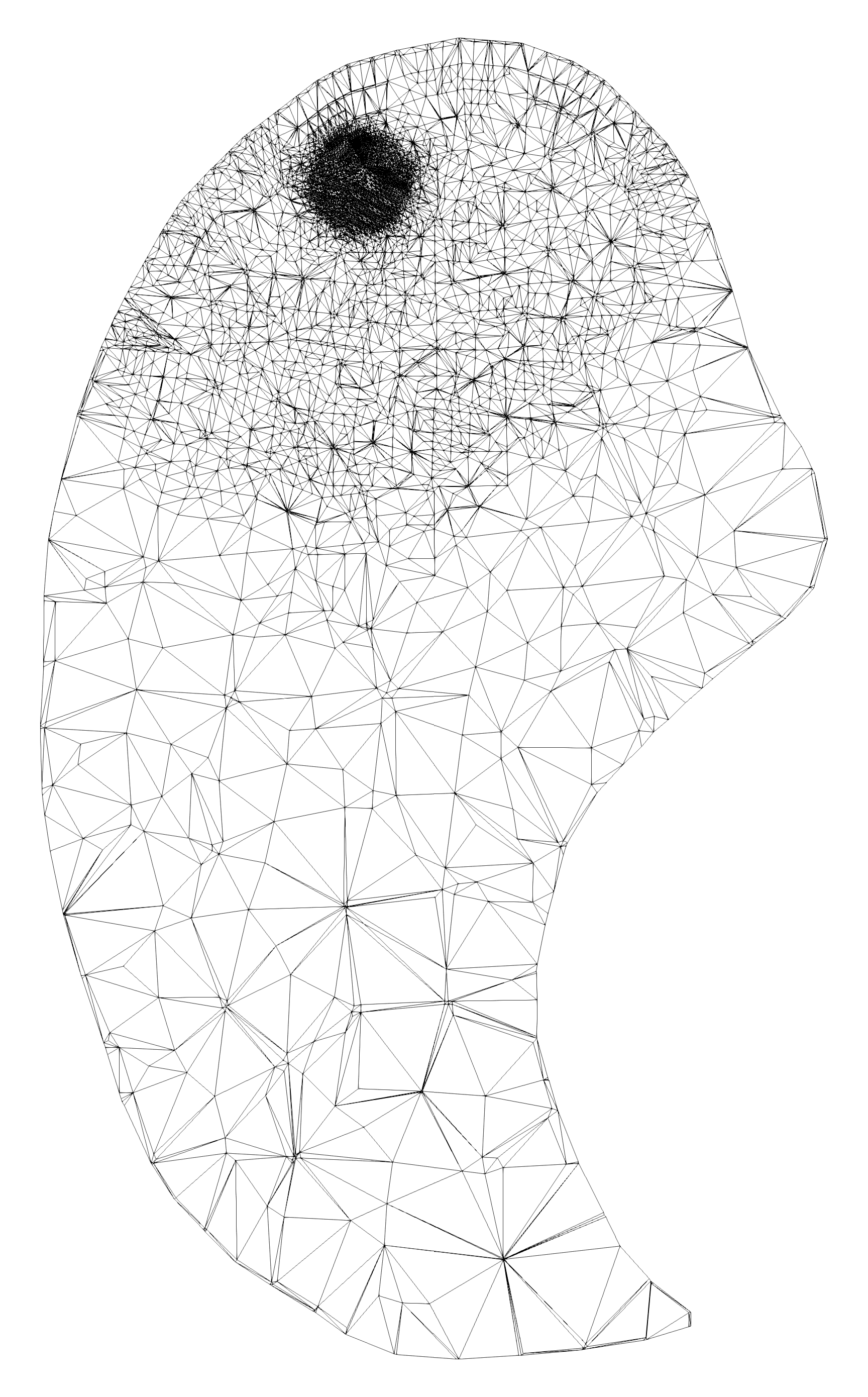}
    \end{minipage}
    \begin{minipage}{.5\textwidth+1.125in}
    \centering (A)
    \end{minipage}%
    \begin{minipage}{.5\textwidth+1.125in}
    \centering (B)
    \end{minipage}
    \parbox{\linewidth}{
    \caption{
    \label{fig:lung-discretization}
    {\bf  Patient specific lung geometry reconstructed from CT images:}
    A:~Lung envelope and 2D-surface of the vasculature.
    B:~2D slice of the 3D mesh near the tumor.
    }
    }
    \end{adjustwidth}
\end{figure}
The vasculature and the lung are depicted in Fig.~\ref{fig:lung-discretization}A.
To keep the computational costs tractable, we apply a local refinement close to the initial tumor position by subdividing the mesh several times, see Fig.~\ref{fig:lung-discretization}B.
This approach is motivated by the potential energy of the Cahn--Hilliard equation, which keeps the central part of the tumor localized to its initial location. 
Similarly, we only require an accurate solution close to the tumor for the nutrient model, which couples with the Cahn--Hilliard equation.  

\section*{Acknowledgements}\label{sec:acc}
The work of A. Wagner and B. Wohlmuth was partially funded by the Deutsche Forschungsgemeinschaft (WO 671/11-1). P. Schlicke was partially funded by the IGSSE/TUM Graduate School.
Marvin Fritz is partially supported by the State of Upper Austria.
The support of J. Tinsley Oden by the U.S. Dept. of Energy, Office of Science, Office of Advanced Scientific Computing Research, Applied Mathematics Program, under Award DE-960009286 is gratefully acknowledged. 

\nolinenumbers

{
  \small	
	\bibliography{literature.bib}
}

\end{document}